\documentclass[11pt]{article}
\usepackage{amsmath,amsthm,amssymb,graphicx}
\usepackage[margin=1in]{geometry}
\usepackage{amsfonts}
\usepackage{hyperref}
\usepackage{tikz}
\usepackage{booktabs}
\usepackage[all]{xy}
\newtheorem{lemma}{Lemma}

\newcommand{\Spec}{\mathrm{Spec}}

\newtheorem{theorem}{Theorem}
\newtheorem{definition}{Definition}
\newtheorem{corollary}{Corollary}

\newtheorem{proposition}{Proposition}
\newtheorem{remark}{Remark}
\newcommand{\HL}{\mathrm{HL}}
\theoremstyle{plain}

\newcommand{\BDC}[1]{#1\otimes K_2}
\DeclareMathOperator{\Gap}{Gap}

\begin{document}

\title{Perfecting the Line Graph}
\author{Hartosh Singh Bal}
\date{}

\maketitle
\renewcommand{\thefootnote}{}
\renewcommand{\thefootnote}{\arabic{footnote}}

\begin{abstract}
We study the doubled edge-stage lift
\[
\HL'_2(G)=L(G\otimes K_2),
\]
the line graph of the canonical bipartite double cover of a graph \(G\). The natural involution
\((u,v)\leftrightarrow(v,u)\) has quotient isomorphic to \(L(G)\), and induces a sector decomposition
\[
\Spec(\HL'_2(G))=\Spec(L(G))\cup\Spec(\mathcal A(G)),
\]
where \(\mathcal A(G)\) is a canonical signed refinement of the line graph. Thus the construction retains
substantial edge-space information through its quotient and antisymmetric sector.

For every input graph, \(\HL'_2(G)\) is perfect, claw-free, and box-perfect. In the regular case we give an
explicit spectral formula, together with quantitative control of the second eigenvalue and spectral gap for
non-bipartite input. Explicit families, including the complete-graph lifts and the Paley lifts, illustrate the
theory; in particular, the Paley lifts furnish an explicit family of regular perfect graphs with controlled
adjacency spectrum and spectral gap.

The construction may be viewed both intrinsically, via ordered-edge adjacency by one-coordinate agreement, and
extrinsically, as the line graph of the canonical double cover. The first viewpoint emphasizes the edge-stage
nature of the lift, while the second supplies the structural proofs used here.
\end{abstract}

\section{Introduction}

The ordinary line graph \(L(G)\) is one of the most classical graph constructions, but it loses information:
distinct graphs may have the same line graph, most famously \(K_3\) and \(K_{1,3}\). The aim of this paper is
to place the line graph inside a larger canonical edge-space object that both recovers \(L(G)\) and carries
additional structure.

The central object is the symmetric lift
\[
\HL'_2(G)=L(G\otimes K_2),
\]
the line graph of the canonical bipartite double cover of \(G\). This may be viewed as a doubled edge-stage
representation of the base graph. It comes equipped with a natural involution
\[
(u,v)\leftrightarrow (v,u),
\]
and the quotient by this involution is exactly the ordinary line graph:
\[
\HL'_2(G)/\iota \cong L(G).
\]
Thus \(L(G)\) appears as the symmetric shadow of a larger doubled construction.

This doubled edge-stage object admits a natural sector decomposition. The symmetric sector recovers \(L(G)\),
while the antisymmetric sector yields a canonical signed graph $\mathcal A(G)$, the antisymmetric line graph.
Accordingly,
\[
\Spec(\HL'_2(G))=\Spec(L(G))\cup\Spec(\mathcal A(G)).
\]

The ingredients entering this construction are individually classical, but their assembly in the canonical doubled edge-stage \( \HL'_2(G)=L(G\otimes K_2) \) leads to a sharper package of consequences: the quotient recovery of \(L(G)\), the symmetric/antisymmetric sector decomposition with antisymmetric companion \(\mathcal A(G)\), explicit spectral formulas in the regular case, and quantitative control of the lifted spectral gap.
The point is not merely to combine classical facts about double covers and line graphs, but to isolate a canonical
doubled edge-space object together with its symmetric quotient and its antisymmetric signed companion. This gives a
new way to construct perfect graphs from arbitrary input graphs without collapsing the line-graph shadow
of the original object. One clean consequence is a Whitney-type refinement: the pair consisting of the
ordinary line graph and the antisymmetric signed sector distinguishes the classical \(K_3/K_{1,3}\)
ambiguity.

A guiding point of view in this paper is that the construction has two equivalent descriptions. One is intrinsic:
vertices are ordered edges, and adjacency is given by a one-coordinate agreement rule. The other is extrinsic:
the same graph is realized as the line graph of the canonical double cover \(G\otimes K_2\). The second
description makes the structural properties and quotient theorem transparent, while the first emphasizes that the
construction is a genuine edge-space refinement rather than merely a derived consequence of the tensor model.

Related spectral questions involving line graphs and Kronecker products have also been studied
by Chauhan and Reddy \cite{chauhan2025}, though in a different order of operations: their focus is on
Kronecker products of line graphs, whereas here the central object is the line graph of the canonical
Kronecker double cover.

We also study the ordered lift \(\HL_2(G)\), but only as a secondary, label-dependent object. Its role here is
to exhibit how a choice of acyclic orientation selects one sheet inside the canonical symmetric lift.

Beyond the basic doubled lift and its sector decomposition, we record explicit regular examples, including the
complete-graph lifts and the Paley lifts, and we note that the symmetric lift can be iterated, producing a
canonical tower within the edge-stage itself.

The antisymmetric sector introduced here is developed intrinsically in separate work \cite{bal2026alg}. In the
present paper, its role is to complete the sector decomposition of the symmetric lift. The Paley examples,
together with the regular spectral formula and its gap consequences, already show that the construction behaves
naturally on structured pseudorandom families and suggest a broader random-regular direction that we do not
pursue here.

\section{The Edge-Stage Lift: Two Viewpoints}

We begin with the edge-stage construction that underlies the whole paper. It admits two equivalent
descriptions. The intrinsic one is in terms of ordered edges with Hamming-type adjacency; the extrinsic one
realizes the same graph as the line graph of the canonical bipartite double cover. The second viewpoint is
particularly useful for proving structural properties such as perfectness, box-perfectness, and the quotient
relationship to the line graph.

\begin{definition}[Ordered Lift \( \mathrm{HL}_2(G) \)]\label{def:ordered-lift}
Let \( G = (V, E) \) be a finite simple graph with an injective labeling \( \phi : V \to \mathbb{N} \). Construct a bipartite graph \( B_{\mathrm{HL}}(G) \) as follows:
\begin{itemize}
    \item The vertex set of \( B_{\mathrm{HL}}(G) \) is \( V' \sqcup V'' \), where \( V' = \{ v' : v \in V \} \) and \( V'' = \{ v'' : v \in V \} \) are two disjoint copies of the vertex set.
    \item For each edge \( \{u,v\} \in E \) with \( \phi(u) < \phi(v) \), include an edge from \( u' \in V' \) to \( v'' \in V'' \).
\end{itemize}
The ordered lift $\HL_2(G)$ is the line graph $L(B_{\HL}(G))$ of this bipartite graph.

\end{definition}

This construction may also be viewed intrinsically. Its vertices are ordered edges \((u,v)\), and
two such vertices are adjacent when they agree in exactly one coordinate, i.e. when they have the
same tail or the same head. Thus the ordered lift is naturally governed by a Hamming-type
adjacency rule on ordered edges.

\begin{lemma}
Each vertex of $\mathrm{HL}_2(G)$ corresponds to a directed edge $(u,v)$ with $\phi(u)<\phi(v)$.
Two such vertices $(u,v)$ and $(x,y)$ are adjacent in $\mathrm{HL}_2(G)$ if and only if the
corresponding bipartite edges $u'\!-\!v''$ and $x'\!-\!y''$ share an endpoint in $B_{\mathrm{HL}}(G)$.
Equivalently, $(u,v)$ is adjacent to $(x,y)$ if and only if either $u=x$ (shared tail) or $v=y$
(shared head). In particular, adjacency is exactly ``agreeing in one coordinate'', i.e. Hamming
distance~$1$ on ordered pairs.
\end{lemma}

\begin{proof}
By definition, $\HL_2(G)$ is the line graph of $B_{\HL}(G)$. Thus two vertices of $\HL_2(G)$
are adjacent if and only if the corresponding edges of $B_{\HL}(G)$ share a common endpoint.
Writing those edges as $u'v''$ and $x'y''$, they share an endpoint if and only if either
$u'=x'$ (equivalently $u=x$) or $v''=y''$ (equivalently $v=y$). This is exactly the
condition that $(u,v)$ and $(x,y)$ agree in one coordinate and differ in the other, i.e.
have Hamming distance one in $V\times V$.
\end{proof}

\begin{definition}[Symmetric Lift \( \mathrm{HL}'_2(G) \)]\label{def:symmetric-lift}
Let \( G = (V, E) \) be a finite simple undirected graph. Construct a bipartite graph \( B_{\mathrm{HL}'}(G) \) as follows:
\begin{itemize}
    \item The vertex set is again \( V' \sqcup V'' \), two disjoint copies of \( V \).
    \item For each edge \( \{u,v\} \in E \), include two edges: one from \( u' \in V' \) to \( v'' \in V'' \), and one from \( v' \in V' \) to \( u'' \in V'' \).
\end{itemize}
The symmetric lift \( \mathrm{HL}'_2(G) \) is the line graph \( L(B_{\mathrm{HL}'}(G)) \) of this bipartite graph.
\end{definition}

\begin{remark}[Canonical double cover / tensor product viewpoint]\label{rem:canonical-double-cover}
The bipartite graph $B_{\mathrm{HL}'}(G)$ is the \emph{canonical (Kronecker) double cover} of $G$,
which can be realized as the categorical/tensor product $G\otimes K_2$
(some authors write $G\times K_2$); see \cite{imrichKlavzar2000}.
Equivalently, $A(G\otimes K_2)=A(G)\otimes A(K_2)$, so if
$\operatorname{Spec}(A(G))=\{\lambda_1,\dots,\lambda_n\}$ then
$\operatorname{Spec}(A(G\otimes K_2))=\{\pm\lambda_1,\dots,\pm\lambda_n\}$.
This classical $\pm$-doubling explains the appearance of the $\pm\lambda_i(G)$ terms in the
regular-spectrum formula for $\mathrm{HL}'_2(G)=L(G\otimes K_2)$, while our main spectral
decomposition Theorem~\ref{thm:spectral-decomposition} is proved intrinsically via the involution
$\iota$ and the associated symmetric/antisymmetric eigenspace splitting.
\end{remark}

The ordered lift \(\HL_2(G)\) is direction-sensitive and depends on the chosen labeling. In the present paper,
it plays only a secondary role, chiefly as a label-dependent induced subgraph of the symmetric lift. It is a
rich structure in its own right and is developed separately; it is mentioned here to highlight the close
connection between the two lifts.

\subsection{Graph-Theoretic Properties of the Ordered Lift}

\begin{theorem}[Line-graph connectivity]\label{thm:linegraph-connectivity}
Let $H$ be a finite graph and let $H^\ast$ be the graph obtained from $H$ by deleting isolated vertices.
Then the line graph $L(H)$ is connected if and only if $H^\ast$ is connected.
\end{theorem}

\begin{proof}
Let $H^\ast$ denote the graph obtained from $H$ by deleting all isolated vertices.

We prove that $L(H)$ is connected if and only if $H^\ast$ is connected.

\smallskip
\noindent\emph{($\Rightarrow$)} Assume $L(H)$ is connected. Let $x,y\in V(H^\ast)$.
Choose edges $e,f\in E(H)$ incident to $x$ and $y$, respectively. Since $L(H)$ is connected,
there is a path $e=e_0,e_1,\dots,e_t=f$ in $L(H)$. By definition of the line graph, consecutive
edges $e_i,e_{i+1}$ share an endpoint in $H$. Thus $e_0,e_1,\dots,e_t$ is an edge-walk in $H$,
and following shared endpoints produces a vertex-walk in $H^\ast$ from $x$ to $y$. Hence
$H^\ast$ is connected.

\smallskip
\noindent\emph{($\Leftarrow$)} Assume $H^\ast$ is connected. Let $e,f\in E(H)$ be any two edges.
Pick endpoints $u\in e$ and $v\in f$; then $u,v\in V(H^\ast)$. Since $H^\ast$ is connected,
there is a vertex-path $u=x_0,x_1,\dots,x_m=v$ in $H^\ast$. For each $i$ choose an edge
$g_i\in E(H)$ with endpoints $\{x_i,x_{i+1}\}$. Then $g_0,g_1,\dots,g_{m-1}$ is a path in
$L(H)$, because $g_i$ and $g_{i+1}$ share the vertex $x_{i+1}$ in $H$. Finally, since $e$ is
incident to $x_0=u$, we have $e\sim g_0$ in $L(H)$, and since $f$ is incident to $x_m=v$, we
have $g_{m-1}\sim f$. Concatenating gives a path from $e$ to $f$ in $L(H)$. Therefore
$L(H)$ is connected.
\end{proof}

\begin{corollary}[Connectivity criterion for the ordered lift]\label{cor:HL2-connectivity}
Let $G$ be a finite graph with injective labeling $\phi$, and let $B_{\HL}(G)$ be the bipartite
graph from Definition~\ref{def:ordered-lift}. Let $B_{\HL}(G)^\ast$ denote the graph obtained from
$B_{\HL}(G)$ by deleting isolated vertices. Then $\HL_2(G)=L(B_{\HL}(G))$ is connected if and only if
$B_{\HL}(G)^\ast$ is connected.
\end{corollary}

\begin{proof}
Apply Theorem~\ref{thm:linegraph-connectivity} with \(H=B_{\HL}(G)\).
\end{proof}

In particular, $\HL_2(G)$ may be disconnected even when $G$ is connected, since
$B_{\HL}(G)$ depends on the labeling $\phi$; see the path example below.

\subsection{Structure of the Symmetric Lift {${HL}'_2(G)$}}

\begin{theorem}\label{thm:symmetric-connectedness}
Let \(G\) be a connected undirected graph with at least one edge. Then:
\begin{enumerate}
  \item If \(G\) is not bipartite, then \(\HL'_2(G)\) is connected.
  \item If \(G\) is bipartite, then \(\HL'_2(G)\) has exactly two connected components.
\end{enumerate}
\end{theorem}
\begin{proof}
Recall that \(\HL'_2(G)=L\!\big(B_{\HL'}(G)\big)\), where \(B_{\HL'}(G)\) is the
bipartite (Kronecker) double cover of \(G\).

\noindent\textbf{Fact A (line-graph connectivity).}
This is exactly Theorem~\ref{thm:linegraph-connectivity}.

\medskip
\noindent\textbf{Fact B (connectivity of the bipartite double cover).}
If \(G\) is connected, then \(B_{\HL'}(G)\) is connected if and only if \(G\) is
not bipartite; if \(G\) is bipartite, then \(B_{\HL'}(G)\) has exactly two
connected components.

\smallskip
We prove Fact B. Write \(V(B_{\HL'}(G))=V'\sqcup V''\), and for each
\(\{u,v\}\in E(G)\) the cover contains the two edges \(u'v''\) and \(v'u''\).

\smallskip
\emph{(Bipartite \(\Rightarrow\) two components).}
Assume \(G\) is bipartite with bipartition \(V(G)=A\sqcup B\).
Then every edge of \(B_{\HL'}(G)\) goes either from \(A'\) to \(B''\) or from
\(B'\) to \(A''\). Hence
\[
H_1 := B_{\HL'}(G)[A'\cup B'']\qquad\text{and}\qquad
H_2 := B_{\HL'}(G)[B'\cup A'']
\]
are spanning subgraphs whose edge sets partition \(E(B_{\HL'}(G))\), and there
are no edges between \(H_1\) and \(H_2\). Since \(G\) is connected and has an
edge, each \(H_i\) contains at least one edge. Moreover, any path in \(G\)
lifts to a path inside each \(H_i\), so each \(H_i\) is connected. Therefore
\(B_{\HL'}(G)\) has exactly two connected components.

\smallskip
\emph{(Non-bipartite \(\Rightarrow\) connected).}
Assume \(G\) is connected and not bipartite. Then \(G\) contains an odd cycle
\(u_0u_1\cdots u_{2t}u_0\). In \(B_{\HL'}(G)\) the lifted edges
\(
u_0'u_1'',\;
u_2'u_1'',\;
u_2'u_3'',\;
u_4'u_3'',\;
\dots,\;
u_{2t}'u_0''
\)
give a path
\[
u_0' \sim u_1'' \sim u_2' \sim u_3'' \sim \cdots \sim u_{2t}' \sim u_0'',
\]
so \(u_0'\) and \(u_0''\) lie in the same connected component.

Now fix any vertex \(x\in V(G)\). Since \(G\) is connected, choose a path
\(u_0=x_0,x_1,\dots,x_m=x\) in \(G\). Lifting it starting at \(u_0'\) yields a
path in \(B_{\HL'}(G)\) from \(u_0'\) to \(x'\) if \(m\) is even and to \(x''\)
if \(m\) is odd. Composing with the established path between \(u_0'\) and
\(u_0''\) toggles the sheet, so both \(x'\) and \(x''\) lie in the same
component as \(u_0'\). Hence \(B_{\HL'}(G)\) is connected. This proves Fact B.

\medskip
Now apply Facts A and B with \(H=B_{\HL'}(G)\). Since \(G\) is connected and has
an edge, the cover \(B_{\HL'}(G)\) has no isolated vertices, so \(H^\ast=H\).
If \(G\) is not bipartite, then \(B_{\HL'}(G)\) is connected by Fact B, hence
\(\HL'_2(G)=L(B_{\HL'}(G))\) is connected by Fact A. If \(G\) is bipartite, then
\(B_{\HL'}(G)\) has exactly two connected components each containing edges, so
\(L(B_{\HL'}(G))\) has exactly two connected components by Fact A.
\end{proof}

\begin{theorem}\label{thm:symmetric-regularity}
If \( G \) is \( d \)-regular, then \( \mathrm{HL}'_2(G) \) is \( (2d - 2) \)-regular.
\end{theorem}

\begin{proof}
Let $G$ be $d$-regular. A vertex of $\HL'_2(G)=L(B_{\HL'}(G))$ corresponds to
an edge of $B_{\HL'}(G)$, say $u'v''$ arising from $\{u,v\}\in E(G)$.
In the line graph, its neighbors are exactly the other edges of $B_{\HL'}(G)$
incident to $u'$ together with the other edges incident to $v''$.

Since $G$ is $d$-regular, $u$ has exactly $d$ neighbors in $G$, so $u'$ is
incident in $B_{\HL'}(G)$ to exactly $d$ edges $u'w''$ (one for each neighbor
$w$ of $u$). Excluding $u'v''$ leaves $d-1$ neighbors in the line graph.
Similarly, $v''$ is incident to exactly $d$ edges $x'v''$ (one for each neighbor
$x$ of $v$), and excluding $u'v''$ leaves another $d-1$ neighbors. These two
sets are disjoint, so
\[
\deg_{\HL'_2(G)}(u'v'')=(d-1)+(d-1)=2d-2.
\]
\end{proof}

\subsubsection{Example: A Path on Three Vertices}
Let $G$ be the path graph on three vertices $a, b, c$, with edges $e_1 = \{a,b\}$ and $e_2 = \{b,c\}$.

\[
\xymatrix@R=1em@C=3em{
a \ar[r]^{e_1} & b \ar[r]^{e_2} & c
}
\]

\paragraph{Case A}
We consider a labeling where $a$ is assigned $0$, $b$ is assigned $1$, and $c$ is assigned $2$.

In the ordered lift $\mathrm{HL}_2(G)$, the double cover consists of vertices $0'$, $1'$, $2'$ and $0''$, $1''$, $2''$. The edges go from the top to the bottom cover following the labeling order in the graph and include $0'1''$ and $1'2''$. Taking the line graph, these correspond to two vertices in the lift with no shared vertex in the bipartite cover, and hence are not adjacent.

\[
\xymatrix@C=4em{
*+[o][F-]{0'1''} & *+[o][F-]{1'2''}
}
\]

\noindent $\mathrm{HL}_2(G)$ is therefore disconnected.

\paragraph{Case B}
We consider a labeling where $a$ is assigned $1$, $b$ is assigned $0$, and $c$ is assigned $2$.

In the ordered lift $\mathrm{HL}_2(G)$, the double cover consists of $0'$, $1'$, $2'$ and $0''$, $1''$, $2''$. The edges are $0'1''$ and $0'2''$, which share the top vertex $0'$, so their corresponding vertices in the lift are adjacent.

\[
\xymatrix@C=4em{
*+[o][F-]{0'1''} \ar@{-}[r] & *+[o][F-]{0'2''}
}
\]

\noindent $\mathrm{HL}_2(G)$ is therefore connected.

\paragraph{Symmetric Lift {$\mathrm{HL}'_2(G)$}}
We construct the symmetric lift $\mathrm{HL}'_2(G)$ explicitly as follows.

\textbf{Step 1: Bipartite double cover.} 

Form two disjoint copies of the vertex set: $V' = \{a', b', c'\}$ and $V'' = \{a'', b'', c''\}$. For each edge in $G$, we add edges between these copies in both directions:

\begin{itemize}
  \item Edge $e_1 = \{a,b\}$ becomes: $(a',b'')$ and $(b',a'')$
  \item Edge $e_2 = \{b,c\}$ becomes: $(b',c'')$ and $(c',b'')$
\end{itemize}

This bipartite double cover is the \emph{canonical/Kronecker double cover} of $G$, $\BDC{G}$: it has vertex set $V(G)\times\{0,1\}$ and edges $(u,0)\sim(v,1)$ iff $u\sim v$ in $G$.

\textbf{Step 2: Forming $\mathrm{HL}'_2(G)$ (line graph).} The vertices of $\mathrm{HL}'_2(G)$ correspond to the edges above. Adjacency is determined by shared vertices in the bipartite double cover:

\begin{itemize}
  \item $(a', b'')$ and $(c', b'')$ share $b''$ — they are adjacent.
  \item $(b', a'')$ and $(b', c'')$ share $b'$ — they are adjacent.
\end{itemize}

There are no other adjacencies, giving two disconnected edges. Thus, the symmetric lift $\mathrm{HL}'_2(G)$ is:

\[
\xymatrix@C=4em{
*+[o][F-]{(a',b'')} \ar@{-}[r] & *+[o][F-]{(c',b'')} & \hspace{4em} & *+[o][F-]{(b',a'')} \ar@{-}[r] & *+[o][F-]{(b',c'')}
}
\]

Note that this result does not depend on how vertices $a$, $b$, or $c$ are labeled. Any labeling yields exactly the same structure for $\mathrm{HL}'_2(G)$. This illustrates the general fact established in Theorem~\ref{thm:symmetric-connectedness} that since the original path graph $G$ is bipartite, its symmetric lift $\mathrm{HL}'_2(G)$ consists of exactly two disconnected components.

\paragraph{Quotient Structure.}
A natural automorphism of order two exists on \( \mathrm{HL}'_2(G) \), exchanging the directed pairs \( (u, v) \leftrightarrow (v, u) \). Taking the quotient by this involution identifies each vertex with its twin, yielding a graph isomorphic to the classical line graph \( L(G) \). The symmetric lift thus serves as a canonical doubled cover of the line graph, retaining perfection and
box-perfectness while encoding \(L(G)\) through this quotient.

\subsection{Recoverability and quotient structure}

We now explain how the base graph structure is recovered from the symmetric lift and its quotient.
The key point is that the symmetric lift carries a canonical involution
\[
(u,v)\leftrightarrow (v,u),
\]
whose quotient recovers the ordinary line graph. Thus the doubled edge-stage lift should be viewed as a canonical
enlargement of \(L(G)\), rather than as a separate unrelated construction.

The ordered lift plays a secondary role here: it should be viewed as a label-dependent sheet selection
inside the canonical symmetric lift. The symmetric lift is the primary object because it is independent
of any orientation choice, admits the involution $(u,v)\leftrightarrow(v,u)$, and naturally decomposes
into symmetric and antisymmetric sectors. The ordered lift survives chiefly as an induced subgraph
obtained by selecting one orientation from each involution pair.

The symmetric lift \(\HL'_2(G)\) canonically encodes the line graph \(L(G)\): the involution interchanging the two
sheets of \(\BDC{G}\) induces a natural \(2\)-cover
\[
\HL'_2(G)\to L(G).
\]
In particular, for connected simple graphs, \(\HL'_2(G)\) determines \(G\) up to the classical Whitney ambiguity:
if \(\HL'_2(G)\cong \HL'_2(H)\), then \(L(G)\cong L(H)\), hence either \(G\cong H\) or \(\{G,H\}=\{K_3,K_{1,3}\}\).
The antisymmetric sector resolves this remaining ambiguity.

\begin{lemma}\label{lem:ordered-lift-acyclic-orientation}
Let $G$ be a graph. The ordered lift $\mathrm{HL}_2(G)$ depends only on the acyclic orientation of
$E(G)$ induced by the labeling $\phi$ (orient each edge from smaller label to larger label).
Consequently, the number of distinct ordered lifts obtainable from $G$ is at most the number of
acyclic orientations of $G$.
\end{lemma}

\begin{proof}
By the definition of $\mathrm{HL}_2(G)$ and the adjacency criterion above, adjacency in
$\mathrm{HL}_2(G)$ is determined entirely by which endpoint of each edge is the tail and which is
the head; that is, by the induced orientation of $E(G)$. Any injective labeling produces an acyclic
orientation, and two labelings that induce the same oriented edges yield the same bipartite graph
$B_{\mathrm{HL}}(G)$ up to relabeling of the copies, hence the same line graph. The final claim
follows.
\end{proof}

\subsection{Structural Properties: Perfectness and Claw-Freeness}

Both the ordered lift \( \mathrm{HL}_2(G) \) and the symmetric lift \( \mathrm{HL}'_2(G) \) are constructed as line graphs of bipartite graphs. This guarantees that they inherit several strong structural properties, regardless of the base graph \( G \) or the choice of labeling.

\begin{theorem}\label{thm:lifts-perfect-clawfree}
For any finite graph \( G \), both \( \mathrm{HL}_2(G) \) and \( \mathrm{HL}'_2(G) \) are perfect,
claw-free, and odd-hole-free graphs.
\end{theorem}

\begin{proof}
Both \( \HL_2(G) \) and \( \HL'_2(G) \) are line graphs of bipartite graphs by construction.
Such graphs are perfect; see \cite{chudnovsky2006spgt,schrijver2003}. They are also odd-hole-free,
since line graphs of bipartite graphs contain no induced odd cycle of length at least five.
Finally, every line graph is claw-free by Beineke’s characterization \cite{beineke1970}.
Therefore both lifts are perfect, claw-free, and odd-hole-free.
\end{proof}

\begin{theorem}[Box-perfectness of the lifts]\label{thm:lifts-box-perfect}
For any finite graph \( G \), both \( \mathrm{HL}_2(G) \) and \( \mathrm{HL}'_2(G) \) are
box-perfect graphs.
\end{theorem}

\begin{proof}
Both lifts are line graphs of bipartite graphs, and this class is box-perfect; see \cite{schrijver2003}.
Therefore both \( \HL_2(G) \) and \( \HL'_2(G) \) are box-perfect.
\end{proof}

Box-perfect graphs form a distinguished subclass of perfect graphs whose stable set polytopes admit particularly strong linear
descriptions under box constraints. The stable set polytope of a graph \( G=(V,E) \) is
\[
\mathrm{STAB}(G)=\mathrm{conv}\bigl\{\mathbf{1}_S\in\mathbb{R}^V : S\subseteq V \text{ is a stable set}\bigr\}.
\]
For line graphs of bipartite graphs, this polyhedral structure is especially well behaved.

\begin{theorem}[Clique Number of Symmetric Lifts]\label{thm:clique-number}
Let \( G \) be a \( d \)-regular graph. Then:
\[
\omega(\mathrm{HL}'_2(G)) = d.
\]
\end{theorem}

\begin{proof}
Let \(G\) be \(d\)-regular. Fix a vertex \(v\in V(G)\). In the bipartite double cover \(B_{\HL'}(G)\),
the vertex \(v'\) is incident to exactly \(d\) edges, namely the edges \(v'u''\) for \(u\sim v\) in \(G\).
These \(d\) edges correspond to \(d\) vertices of the line graph \(\HL'_2(G)\), and because they are all
incident to the common vertex \(v'\), they form a clique of size \(d\).

To see that no larger clique exists, recall that \(\HL'_2(G)\) is the line graph of the bipartite graph
\(B_{\HL'}(G)\). In a bipartite graph, every clique in the line graph comes from a family of edges all
incident to a single vertex of the underlying graph, since a bipartite graph contains no triangles.
Thus the largest clique in \(\HL'_2(G)\) has size equal to the maximum degree of \(B_{\HL'}(G)\), which is
\(d\). Therefore \(\omega(\HL'_2(G))=d\).
\end{proof}

The perfection of the symmetric lift yields:

\begin{corollary}[Chromatic Number of Symmetric Lifts]
Let \( G \) be a \( d \)-regular graph. Then the symmetric lift \( \mathrm{HL}'_2(G) \) satisfies
\[
\chi(\mathrm{HL}'_2(G)) = d.
\]
\end{corollary}

\begin{proof}
By Theorem~\ref{thm:clique-number}, \(\omega(\HL'_2(G))=d\). By Theorem~\ref{thm:lifts-perfect-clawfree},
\(\HL'_2(G)\) is perfect. Hence for the perfect graph \(\HL'_2(G)\), the chromatic number equals the
clique number:
\[
\chi(\HL'_2(G))=\omega(\HL'_2(G))=d.
\]
\end{proof}

The symmetric lift is the canonical object in the present paper; the ordered lift should be viewed as a
label-dependent sheet selection inside it.

\begin{proposition}\label{prop:HLprime2-is-linegraph}
For every graph $G$, the symmetric lift $\HL'_2(G)$ is the line graph of the
bipartite double cover $B_{\HL'}(G)$. Equivalently, $\HL'_2(G)$ has vertex set
\[
V(\HL'_2(G))=\{(u,v):\{u,v\}\in E(G)\},
\]
and two vertices $(u,v)$ and $(x,y)$ are adjacent in $\HL'_2(G)$ if and only if
either $u=x$ or $v=y$ when interpreted as bipartite edges $u'v''$ and $x'y''$.
\end{proposition}

\begin{proof}
By definition, $B_{\HL'}(G)$ has bipartition $V'\sqcup V''$ and for each
$\{u,v\}\in E(G)$ contains the two edges $u'v''$ and $v'u''$. The line graph
$L(B_{\HL'}(G))$ therefore has a vertex for each such bipartite edge, which we
identify with the ordered pair $(u,v)$ corresponding to $u'v''$.
Two vertices of the line graph are adjacent exactly when the corresponding
bipartite edges share an endpoint, i.e.\ share $u'$ (so $u=x$) or share $v''$
(so $v=y$). This is precisely the stated adjacency rule, and yields
$\HL'_2(G)=L(B_{\HL'}(G))$.
\end{proof}

\subsection{The symmetric lift and the line graph}

\begin{proposition}\label{prop:HL2-induced-in-HLprime2}
Let $G$ be a graph with an injective labeling $\phi:V(G)\to\mathbb{N}$, and let
\[
S_\phi \;:=\; \{(u,v)\in V(\HL'_2(G)):\ \{u,v\}\in E(G)\ \text{and}\ \phi(u)<\phi(v)\}.
\]
Then the ordered lift $\HL_2(G)$ is isomorphic to the induced subgraph
$\HL'_2(G)[S_\phi]$.
\end{proposition}

\begin{proof}
By Definition~\ref{def:ordered-lift}, the vertices of $\HL_2(G)$ are exactly the oriented edges
$(u,v)$ for which $\{u,v\}\in E(G)$ and $\phi(u)<\phi(v)$. These are precisely the vertices in
the set \(S_\phi\subseteq V(\HL'_2(G))\).

Now let $(u,v),(x,y)\in S_\phi$. In the symmetric lift, they are adjacent exactly when they agree
in one coordinate, equivalently when the corresponding bipartite edges $u'v''$ and $x'y''$ share
an endpoint. But this is exactly the adjacency rule defining the ordered lift. Therefore the
adjacency relation on the induced subgraph $\HL'_2(G)[S_\phi]$ is exactly the adjacency relation
of $\HL_2(G)$.

Hence $\HL_2(G)\cong \HL'_2(G)[S_\phi]$.
\end{proof}

Thus the symmetric lift is the canonical object, while the ordered lift is a label-dependent sheet selection
inside it. The examples \(K_3\) and \(K_{1,3}\) illustrate the point: although their line graphs are
isomorphic, the antisymmetric sector distinguishes them.

\section{Spectral Decomposition of the Symmetric Lift}

We begin with the structural fact underlying the whole theory: the symmetric lift admits a canonical involution
\((u,v)\leftrightarrow (v,u)\), and the quotient by this involution is the ordinary line graph.
This identifies \(L(G)\) as the symmetric shadow of the doubled edge-stage construction.

\begin{theorem}\label{thm:quotient-is-linegraph}
Let $\iota$ be the involution on $\HL'_2(G)$ given by $\iota(u,v)=(v,u)$.
Let $\HL'_2(G)/\iota$ be the quotient graph whose vertices are the $\iota$-orbits.
Then $\HL'_2(G)/\iota \cong L(G)$.
\end{theorem}

\begin{proof}
A vertex of $\HL'_2(G)$ is an ordered pair $(u,v)$ with $\{u,v\}\in E(G)$, and
$\iota$ identifies $(u,v)$ with $(v,u)$. Thus the vertices of the quotient are in
bijection with unordered edges $\{u,v\}\in E(G)$, i.e.\ with $V(L(G))$.

It remains to check adjacency. Two distinct orbits
$[\!(u,v)\!]$ and $[\!(x,y)\!]$ are adjacent in the quotient if and only if
there exist representatives from each orbit that are adjacent in $\HL'_2(G)$.
But adjacency in $\HL'_2(G)=L(B_{\HL'}(G))$ means the corresponding bipartite
edges share an endpoint, i.e.\ for representatives $(a,b)$ and $(c,d)$ we have
$a=c$ or $b=d$.

Now suppose the underlying unordered edges $\{u,v\}$ and $\{x,y\}$ share a
vertex in $G$, say $v=x$ with $u\neq y$. Consider the representatives $(u,v)$
and $(y,v)=(v,y)$ under $\iota$. In $B_{\HL'}(G)$ these correspond to bipartite
edges $u'v''$ and $y'v''$, which share the endpoint $v''$, hence are adjacent in
$\HL'_2(G)$. Therefore the orbits are adjacent in the quotient.

Conversely, if two orbits are adjacent in the quotient, pick adjacent
representatives $(a,b)$ and $(c,d)$. If $a=c$, then the base edges
$\{a,b\}$ and $\{a,d\}$ share the vertex $a$; if $b=d$, then $\{a,b\}$ and
$\{c,b\}$ share the vertex $b$. Hence the corresponding vertices of $L(G)$ are
adjacent.

Thus orbit adjacency is exactly “share an endpoint in $G$”, which is the
adjacency relation of $L(G)$.
\end{proof}

This quotient relationship is the structural bridge between the doubled lift and the ordinary line graph.
The next step is to identify the full sector decomposition induced by the involution.

\subsection{Symmetric and antisymmetric spectral sectors}

\begin{theorem}[Spectral decomposition via symmetric and antisymmetric sectors]\label{thm:spectral-decomposition}
Let \( G \) be a finite simple undirected graph. Then the adjacency spectrum of the symmetric lift
\( \mathrm{HL}'_2(G) \) is the multiset union of the adjacency spectrum of the ordinary line graph
\( L(G) \) and the spectrum of the antisymmetric sector \(\mathcal A(G)\), represented in any oriented-edge
basis as a signed graph well-defined up to switching. That is,
\[
\operatorname{Spec}(\mathrm{HL}'_2(G)) = \operatorname{Spec}(L(G)) \cup \operatorname{Spec}(\mathcal A(G)),
\]
with multiplicities.
\end{theorem}

\begin{proof}
Let \(V=\mathbb R^{V(\HL'_2(G))}\), with basis vectors indexed by ordered edges \((u,v)\) of \(G\).
The involution
\[
\iota:(u,v)\longmapsto (v,u)
\]
acts linearly on \(V\) and commutes with the adjacency operator of \(\HL'_2(G)\), since the adjacency rule
“agree in one coordinate” is preserved by reversing both ordered edges. Hence \(V\) decomposes as a direct sum
of \(\iota\)-invariant subspaces
\[
V=V_+\oplus V_-,
\]
where
\[
V_+=\{f:\ f(u,v)=f(v,u)\},
\qquad
V_-=\{f:\ f(u,v)=-f(v,u)\}.
\]

Fix once and for all an orientation of each undirected edge \(e=\{u,v\}\), writing it as \(u\to v\). For
each such oriented edge \(e\), define
\[
s_e=\delta_{(u,v)}+\delta_{(v,u)}\in V_+,
\qquad
a_e=\delta_{(u,v)}-\delta_{(v,u)}\in V_-.
\]
Then \(\{s_e\}_{e\in E(G)}\) is a basis of \(V_+\), and \(\{a_e\}_{e\in E(G)}\) is a basis of \(V_-\).

On the symmetric basis, two basis vectors \(s_e\) and \(s_f\) interact exactly when the underlying edges
\(e\) and \(f\) share a vertex in \(G\). Thus the matrix of the restricted operator on \(V_+\) is precisely
the adjacency matrix of the ordinary line graph \(L(G)\).

On the antisymmetric basis, the matrix entry between \(a_e\) and \(a_f\) is:
\[
+1 \text{ if } e,f \text{ meet at the same end (tail--tail or head--head),}
\]
\[
-1 \text{ if } e,f \text{ meet at opposite ends (tail--head or head--tail),}
\]
and \(0\) otherwise. This is exactly the signed overlap matrix defining the antisymmetric sector
\(\mathcal A(G)\).

If one changes the chosen orientation of a single edge \(e\), then \(a_e\) is replaced by \(-a_e\), while the
other basis vectors remain unchanged. Thus the antisymmetric block changes by conjugation with a diagonal
\(\pm1\)-matrix, so the resulting signed graph is well-defined up to switching.

Therefore the adjacency operator of \(\HL'_2(G)\) is block diagonal with one block equal to the adjacency
matrix of \(L(G)\) and the other equal to a switching representative of \(\mathcal A(G)\). Consequently,
\[
\Spec(\HL'_2(G))
=
\Spec(L(G))
\cup
\Spec(\mathcal A(G)),
\]
with multiplicities.
\end{proof}
\begin{lemma}[Switching well-definedness of $\mathcal A(G)$]\label{lem:switching-welldefined}
Let $G$ be a finite simple graph and fix any choice of orientation of its edges, so that each
undirected edge $\{u,v\}$ is represented by exactly one oriented edge $e=(u,v)$.
Let $M$ be the signed adjacency matrix on the set of oriented edges defined by the overlap rules:
$M_{e,f}=+1$ for same-end overlaps (tail--tail or head--head), $M_{e,f}=-1$ for crossed overlaps
(tail--head or head--tail), and $M_{e,f}=0$ otherwise.
If the orientation is changed by reorienting an arbitrary subset of edges, producing a new matrix
$M'$, then there exists a diagonal matrix $D$ with $\pm1$ entries such that
\[
M' = D M D.
\]
Equivalently, the signed graph defined by $M$ is well-defined up to switching.
\end{lemma}

\begin{proof}
Reorienting a single edge $e$ replaces its oriented representative by the opposite orientation.
Let $s(e)\in\{\pm1\}$ be $-1$ if $e$ is reoriented and $+1$ otherwise, and let
$D=\mathrm{diag}(s(e))$.
We claim that for every pair of oriented edges $e,f$,
\[
M'_{e,f} = s(e)\,s(f)\,M_{e,f}.
\]
This implies $M'=DM D$.

If neither $e$ nor $f$ is reoriented, then $M'_{e,f}=M_{e,f}$ and the formula holds.
If exactly one of $e,f$ is reoriented, say $e$, then the overlap type between $e$ and $f$ flips
from “same-end” to “crossed” or vice versa, because reversing $e$ swaps its head and tail.
Thus the sign of the interaction changes: $M'_{e,f}=-M_{e,f}$, which equals $s(e)s(f)M_{e,f}$
since $s(e)=-1$ and $s(f)=+1$.
If both $e$ and $f$ are reoriented, then both heads/tails swap and the overlap type is preserved,
so $M'_{e,f}=M_{e,f}=s(e)s(f)M_{e,f}$.
Therefore $M'=DMD$, proving switching equivalence.
\end{proof}

\subsection{The Regular Case}

\begin{theorem}[Explicit Spectrum of the Symmetric Lift for Regular Graphs]\label{thm:regular-spectrum}
Let \( G \) be a finite, simple, \( d \)-regular graph on \(n\) vertices, with adjacency eigenvalues
\( \lambda_1 = d \ge \lambda_2 \ge \cdots \ge \lambda_n \). Then
\[
\operatorname{Spec}(\mathrm{HL}'_2(G))
=
\{\, d - 2 + \lambda_i \mid i=1,\dots,n \,\}
\cup
\{\, d - 2 - \lambda_i \mid i=1,\dots,n \,\}
\cup
\{-2\}^{\times\, n(d-2)}.
\]
In particular, this accounts for all \(nd=2|E(G)|\) eigenvalues of \(\HL'_2(G)\).
\end{theorem}

\begin{proof}
By Theorem~\ref{thm:spectral-decomposition}, the adjacency spectrum of $\HL'_2(G)$ is the multiset union
of the spectra of $L(G)$ and $\mathcal A(G)$.

For a $d$-regular graph $G$ on $n$ vertices, with $m=|E(G)|=nd/2$, the ordinary line graph has spectrum
\[
\Spec(L(G))
=
\{\, \lambda_i+d-2 : 1\le i\le n \,\}
\cup
\{-2\}^{\times (m-n)},
\]
a standard formula for line graphs of regular graphs; see, for example, \cite{godsil2001algebraic}.

To compute the antisymmetric-sector spectrum, fix an orientation of the edges of $G$ and let
$D\in\mathbb R^{n\times m}$ be the corresponding oriented incidence matrix. For distinct edges $e,f$,
the $(e,f)$-entry of $D^\top D$ is
\[
(D^\top D)_{e,f}=
\begin{cases}
+1,& \text{if } e,f \text{ meet with the same local orientation (tail--tail or head--head)},\\
-1,& \text{if } e,f \text{ meet with opposite local orientation (tail--head or head--tail)},\\
0,& \text{otherwise},
\end{cases}
\]
while the diagonal entries are \(2\). Hence
\[
D^\top D-2I
\]
is exactly the signed adjacency matrix of \(\mathcal A(G)\) in the chosen oriented-edge basis. Therefore
\[
\Spec(\mathcal A(G))=\Spec(D^\top D-2I).
\]

Now \(DD^\top\) is the combinatorial Laplacian of \(G\), so
\[
DD^\top = dI-A(G),
\]
where \(A(G)\) is the adjacency matrix of \(G\). The nonzero eigenvalues of \(D^\top D\) and \(DD^\top\)
coincide, so
\[
\Spec(D^\top D)
=
\{\, d-\lambda_i : 1\le i\le n \,\}
\cup
\{0\}^{\times (m-n)}.
\]
Subtracting \(2I\) gives
\[
\Spec(\mathcal A(G))
=
\{\, d-2-\lambda_i : 1\le i\le n \,\}
\cup
\{-2\}^{\times (m-n)}.
\]

Taking the multiset union of the two sector spectra gives
\[
\Spec(\HL'_2(G))
=
\{\, d-2+\lambda_i : 1\le i\le n \,\}
\cup
\{\, d-2-\lambda_i : 1\le i\le n \,\}
\cup
\{-2\}^{\times 2(m-n)}.
\]
Since \(2(m-n)=nd-2n=n(d-2)\), this is exactly
\[
\{\, d - 2 + \lambda_i \mid i=1,\dots,n \,\}
\cup
\{\, d - 2 - \lambda_i \mid i=1,\dots,n \,\}
\cup
\{-2\}^{\times\, n(d-2)}.
\]
This accounts for all \(2m=nd\) eigenvalues of \(\HL'_2(G)\).
\end{proof}

\begin{theorem}[Second eigenvalue and spectral gap of the symmetric lift]\label{thm:lift-second-eigenvalue}
Let \(G\) be a connected, non-bipartite, \(d\)-regular graph, with adjacency eigenvalues
\[
\lambda_1=d>\lambda_2\ge \cdots \ge \lambda_n.
\]
Then \(\HL'_2(G)\) is a connected \((2d-2)\)-regular graph, and its second-largest adjacency eigenvalue is
\[
\lambda_2\!\big(\HL'_2(G)\big)
=
\max\{\,d-2+\lambda_2,\ d-2-\lambda_n\,\}.
\]
Equivalently, its adjacency spectral gap is
\[
\Gap\!\big(\HL'_2(G)\big)
=
(2d-2)-\lambda_2\!\big(\HL'_2(G)\big)
=
\min\{\,d-\lambda_2,\ d+\lambda_n\,\}.
\]
\end{theorem}

\begin{proof}
By Theorem~\ref{thm:regular-spectrum}, the adjacency spectrum of \(\HL'_2(G)\) consists of
\[
d-2+\lambda_1=2d-2,
\]
the values \(d-2+\lambda_i\) for \(i\ge 2\), the values \(d-2-\lambda_i\) for all \(i\), and the flat band
\(-2\).

Since \(G\) is connected and non-bipartite, Theorem~\ref{thm:symmetric-connectedness} implies that
\(\HL'_2(G)\) is connected, hence its top eigenvalue \(2d-2\) has multiplicity one. Therefore the
second-largest eigenvalue is the largest element of the remaining spectrum. Among the values
\(d-2+\lambda_i\), the largest is \(d-2+\lambda_2\). Among the values \(d-2-\lambda_i\), the largest is
\(d-2-\lambda_n\). Since \(G\) is non-bipartite, we have \(\lambda_n>-d\), so
\(d-2-\lambda_n>-2\). Hence the flat band at \(-2\) does not contribute the second-largest eigenvalue.
\[
\lambda_2\!\big(\HL'_2(G)\big)
=
\max\{\,d-2+\lambda_2,\ d-2-\lambda_n\,\}.
\]
Subtracting from the degree \(2d-2\) gives
\[
\Gap\!\big(\HL'_2(G)\big)
=
\min\{\,d-\lambda_2,\ d+\lambda_n\,\}.
\]
\end{proof}

\begin{corollary}[Ramanujan-type input]\label{cor:lift-ramanujan-type}
Let \(G\) be a connected, non-bipartite, \(d\)-regular graph such that
\[
|\lambda_i|\le \theta
\qquad\text{for all } i\ge 2.
\]
Then
\[
\lambda_2\!\big(\HL'_2(G)\big)\le d-2+\theta,
\qquad
\Gap\!\big(\HL'_2(G)\big)\ge d-\theta.
\]
Equivalently, the normalized adjacency spectral gap satisfies
\[
\frac{\Gap(\HL'_2(G))}{2d-2}\ge \frac{d-\theta}{2d-2}.
\]

In particular, if \(G\) is Ramanujan and non-bipartite, so that \(\theta=2\sqrt{d-1}\), then
\[
\lambda_2\!\big(\HL'_2(G)\big)\le d-2+2\sqrt{d-1},
\qquad
\Gap\!\big(\HL'_2(G)\big)\ge d-2\sqrt{d-1}.
\]
\end{corollary}

\begin{proof}
Under the hypothesis \(|\lambda_i|\le \theta\) for \(i\ge 2\), we have
\[
\lambda_2\le \theta
\qquad\text{and}\qquad
-\lambda_n\le \theta.
\]
Applying Theorem~\ref{thm:lift-second-eigenvalue} gives
\[
\lambda_2\!\big(\HL'_2(G)\big)
=
\max\{\,d-2+\lambda_2,\ d-2-\lambda_n\,\}
\le d-2+\theta.
\]
Subtracting from \(2d-2\) yields
\[
\Gap\!\big(\HL'_2(G)\big)\ge d-\theta.
\]
The normalized estimate is immediate, and the Ramanujan case follows by substituting
\(\theta=2\sqrt{d-1}\).
\end{proof}

\section{Detailed structure of the antisymmetric sector}\label{sec:antisym}

The antisymmetric sector of the symmetric lift gives rise to a canonical signed graph associated to \(G\).
While the symmetric sector recovers the ordinary line graph \(L(G)\), the antisymmetric sector records whether
two incident edges meet with compatible or crossed orientation.

This signed graph, denoted \(\mathcal A(G)\), is switching-well-defined and should be viewed as the antisymmetric
companion to the line graph. In the present paper we record only its most immediate structural consequences;
its intrinsic theory is developed separately.

Let \(G=(V,E)\) be a finite simple undirected graph with a total ordering on \(V\). For each edge
\(\{u,v\}\in E\) with \(u<v\), define the oriented edge \((u,v)\), and let \(E^\rightarrow\) denote the
set of such oriented edges. Define a symmetric matrix \(M\in\mathbb{R}^{|E|\times |E|}\), indexed by
\(E^\rightarrow\), by
\begin{itemize}
    \item \(M_{(u,v),(x,y)}=+1\) for same-end overlaps (shared tail or shared head),
    \item \(M_{(u,v),(x,y)}=-1\) for crossed overlaps,
    \item \(M_{(u,v),(x,y)}=0\) otherwise.
\end{itemize}
This is the signed adjacency matrix of \(\mathcal A(G)\).

\begin{theorem}[Rule-Based Matrix for the Antisymmetric Subspace]
Let $G$ be a finite simple undirected graph with a fixed total ordering on $V(G)$. Let $M$ be the matrix defined above. Then $M$ represents the action of the adjacency operator $A$ of $\mathrm{HL}'_2(G)$ restricted to the antisymmetric subspace
\[
\mathcal{A} = \left\{ f: V(\mathrm{HL}'_2(G)) \to \mathbb{R} \mid f(u,v) = -f(v,u) \right\}.
\]
\end{theorem}

\begin{remark}
The matrix identification also follows abstractly from Theorem~\ref{thm:spectral-decomposition} via the
antisymmetric basis \(\{a_e\}\). We include the direct computation below for explicitness.
\end{remark}

\begin{proof}
Let \(G=(V,E)\) be a finite undirected graph with a fixed total ordering on \(V\), inducing oriented edges
\(\{(u,v)\mid u<v,\ \{u,v\}\in E\}\), which index the rows and columns of \(M\). The symmetric lift
\(\HL'_2(G)\) is the line graph of the bipartite double cover \(B\) of \(G\), where \(B\) has vertex sets
\(V'\sqcup V''\) and, for each \(\{u,v\}\in E\), edges \(u'v''\) and \(v'u''\). Vertices of \(\HL'_2(G)\)
are these directed edges, with adjacency when they share a vertex in \(B\), equivalently when they have the
same tail or the same head.

The antisymmetric subspace under \(\sigma:(u,v)\mapsto (v,u)\) consists of functions \(f\) on
\(V(\HL'_2(G))\) satisfying \(f((u,v))=-f((v,u))\). A basis is \(\{e_k\}\), where for the \(k\)-th oriented
edge \((u,v)\) with \(u<v\), one sets \(e_k=\delta_{(u,v)}-\delta_{(v,u)}\). The matrix representation of the
adjacency operator \(A\) of \(\HL'_2(G)\) on this subspace has entries
\[
M_{i,j}=\frac{e_i^{T}Ae_j}{\|e_i\|^2},
\]
with \(\|e_i\|^2=2\). We compute \(e_i^TAe_j\).

For the $i$-th edge $(u,v)$ ($u < v$), $e_i = \delta_{(u,v)} - \delta_{(v,u)}$; for the $j$-th edge $(x,y)$ ($x < y$), $e_j = \delta_{(x,y)} - \delta_{(y,x)}$. Thus,
\[
e_i^T A e_j = [\delta_{(u,v)} - \delta_{(v,u)}]^T A [\delta_{(x,y)} - \delta_{(y,x)}] = \delta_{(u,v)}^T A \delta_{(x,y)} - \delta_{(u,v)}^T A \delta_{(y,x)} - \delta_{(v,u)}^T A \delta_{(x,y)} + \delta_{(v,u)}^T A \delta_{(y,x)}.
\]
Each term $\delta_p^T A \delta_q$ is 1 if directed edges $p$ and $q$ are adjacent in $\mathrm{HL}_2'(G)$, else 0. Adjacency occurs if they share a tail (first coordinates equal, seconds differ) or head (seconds equal, firsts differ). Assume $\{u,v\}$ and $\{x,y\}$ share exactly one vertex (else $M_{i,j} = 0$, matching the rules for no share).

\textbf{Case 1: Shared tail ($u = x$, $v \neq y$)}.  
- $\delta_{(u,v)}^T A \delta_{(x,y)} = \delta_{(u,v)}^T A \delta_{(u,y)} = 1$ (tails $u = u$, heads $v \neq y$ in $B$).
- $\delta_{(u,v)}^T A \delta_{(y,u)} = 0$ (tail $u$ vs head $u$, head $v$ vs tail $y$, no adjacency unless $v = y$, but $v \neq y$).
- $\delta_{(v,u)}^T A \delta_{(x,y)} = \delta_{(v,u)}^T A \delta_{(u,y)} = 0$ (head $u$ vs head $y$, tail $v$ vs tail $u$, no match).
- $\delta_{(v,u)}^T A \delta_{(y,x)} = \delta_{(v,u)}^T A \delta_{(y,u)} = 1$ (heads $u = u$, tails $v \neq y$).
- Total: $1 - 0 - 0 + 1 = 2$. Normalized: $2 / \|e_i\|^2 = 2/2 = 1$, matching the rule (+1 for shared tail).

\textbf{Case 2: Shared head ($v = y$, $u \neq x$)}.  
- $\delta_{(u,v)}^T A \delta_{(x,y)} = \delta_{(u,v)}^T A \delta_{(x,v)} = 1$ (heads $v = v$, tails $u \neq x$).
- $\delta_{(u,v)}^T A \delta_{(y,x)} = \delta_{(u,v)}^T A \delta_{(v,x)} = 0$ (tail $u$ vs head $x$, head $v$ vs tail $v$, no match).
- $\delta_{(v,u)}^T A \delta_{(x,y)} = \delta_{(v,u)}^T A \delta_{(x,v)} = 0$ (head $u$ vs head $v$, tail $v$ vs tail $x$, no match).
- $\delta_{(v,u)}^T A \delta_{(y,x)} = \delta_{(v,u)}^T A \delta_{(v,x)} = 1$ (tails $v = v$, heads $u \neq x$).
- Total: $1 - 0 - 0 + 1 = 2$. Normalized: $1$, matching the rule (+1 for shared head).

\textbf{Case 3: Crossed share ($u = y$, $v \neq x$)}.  
- $\delta_{(u,v)}^T A \delta_{(x,y)} = \delta_{(u,v)}^T A \delta_{(x,u)} = 0$ (tail $u$ vs tail $x$, head $v$ vs head $u$, no match unless $v = u$, impossible).
- $\delta_{(u,v)}^T A \delta_{(y,x)} = \delta_{(u,v)}^T A \delta_{(u,x)} = 1$ (tails $u = u$, heads $v \neq x$).
- $\delta_{(v,u)}^T A \delta_{(x,y)} = \delta_{(v,u)}^T A \delta_{(x,u)} = 1$ (heads $u = u$, tails $v \neq x$).
- $\delta_{(v,u)}^T A \delta_{(y,x)} = \delta_{(v,u)}^T A \delta_{(u,x)} = 0$ (head $u$ vs head $x$, tail $v$ vs tail $u$, no match).
- Total: $0 - 1 - 1 + 0 = -2$. Normalized: $-2/2 = -1$, matching the rule (-1 for crossed share).

\textbf{Case 4: Crossed share ($v = x$, $u \neq y$)}. Symmetric to Case 3, yields $-1$.

\textbf{Case 5: No share or same edge}. All terms are 0 (no adjacencies), so $M_{i,j} = 0$, matching the rule.

Since the matrix entries of $M$ agree with the action of the adjacency operator on every basis vector
of the antisymmetric subspace, $M$ is exactly the matrix of the adjacency operator restricted to that
subspace. In particular, the antisymmetric sector is canonically encoded by the signed graph
$\mathcal A(G)$ up to switching equivalence.
\end{proof}

\subsection{Relation to the ordinary line graph}\label{subsec:connection-line-graph}

The ordinary line graph \(L(G)\) and the antisymmetric signed graph \(\mathcal A(G)\) should be viewed as the
two sector-level shadows of the symmetric lift. The former captures the symmetric overlap structure of incident
edges, while the latter captures the antisymmetric overlap structure.

Taken together, these two sector objects recover more information than either one alone. In particular, they
refine the usual Whitney theory by resolving the \(K_3/K_{1,3}\) ambiguity.

The matrix \(M\) is symmetric with all eigenvalues bounded below by \(-2\). Indeed,
\[
M+2I=(T-H)^\top(T-H),
\]
where \(T\) and \(H\) are the tail and head incidence matrices for oriented edges of \(G\). Thus \(M+2I\)
is a Gram matrix and is therefore positive semidefinite, so \(M\succeq -2I\).
This places the antisymmetric sector in the classical \(-2\)-bounded line-graph/Gram-matrix setting,
while the new feature here is the sign pattern arising from the symmetric/antisymmetric decomposition.

\begin{theorem}[Whitney refinement via the antisymmetric sector]\label{thm:Whitney-refined}
Let \(G\) and \(H\) be connected simple graphs. If
\[
L(G)\cong L(H)
\quad\text{and}\quad
[\mathcal A(G)]=[\mathcal A(H)],
\]
then \(G\cong H\).
\end{theorem}

\begin{proof}
Suppose first that $L(G)\cong L(H)$. By Whitney's classical theorem for connected simple graphs, either
$G\cong H$, or $\{G,H\}=\{K_3,K_{1,3}\}$.

Thus the only remaining ambiguity is the exceptional pair $K_3$ and $K_{1,3}$.
For $K_3$, the antisymmetric signed graph $\mathcal A(K_3)$ is a nontrivial signed triangle; for
$K_{1,3}$, the antisymmetric signed graph consists of three isolated vertices. In particular,
these two signed graphs are not switching-equivalent.

Hence if both $L(G)\cong L(H)$ and $[\mathcal A(G)]=[\mathcal A(H)]$, then the Whitney exceptional case
cannot occur, so necessarily $G\cong H$.
\end{proof}

\section{Regular Graphs and Explicit Examples}

\subsection{Symmetric Lifts of Complete Graphs: A Perfect Refinement of Johnson Graphs}

Let \(K_n\) denote the complete graph on \(n\) vertices. Its symmetric lift \(\HL'_2(K_n)\)
provides a clean explicit family illustrating the general theory developed above. In particular,
it gives a regular perfect graph with a transparent quotient map to \(J(n,2)\) and an explicit
four-eigenvalue spectrum.

\begin{theorem}\label{thm:Kn-lift}
For every $n \ge 4$, the symmetric lift $\mathrm{HL}'_2(K_n)$ is a connected, regular, vertex-transitive,
and edge-transitive graph with $n(n-1)$ vertices and degree $2(n-2)$. It is perfect and box-perfect,
and admits a canonical $2$-to-$1$ projection onto the Johnson graph $J(n,2)$. Its adjacency spectrum
has exactly four distinct eigenvalues,
\[
\mathrm{Spec}(\mathrm{HL}'_2(K_n)) = \{ 2(n-2),\; n-2,\; n-4,\; -2 \},
\]
with multiplicities determined by Theorem~\ref{thm:regular-spectrum}.
\end{theorem}

\begin{proof}
By construction, each vertex of $\mathrm{HL}'_2(K_n)$ corresponds to an ordered pair $(u,v)$ with $u \ne v$ in $K_n$, equivalently an oriented edge $u\to v$. Hence the total number of vertices is $n(n-1)$. The bipartite graph $B_{\mathrm{HL}'}(K_n)$ underlying the lift is formed by taking two disjoint copies of the vertex set of $K_n$ and connecting $u' \in V'$ to $v'' \in V''$ and $v' \in V'$ to $u'' \in V''$ for each edge $\{u,v\}$ in $K_n$. This bipartite graph is connected, regular of degree $n-1$, and its line graph $\mathrm{HL}'_2(K_n)$ is therefore regular of degree $2(n-2)$.

The map $(u,v)\mapsto \{u,v\}$ defines a canonical $2$-to-$1$ graph homomorphism
$\HL'_2(K_n)\to J(n,2)$, identifying the two orientations of each unordered pair.

Since $K_n$ is vertex- and edge-transitive under the symmetric group $S_n$, these symmetries lift naturally to
$\mathrm{HL}'_2(K_n)$, which inherits both vertex- and edge-transitivity. Furthermore, because
$\mathrm{HL}'_2(K_n)$ is the line graph of a bipartite graph, it is perfect and box-perfect by
Theorems~\ref{thm:lifts-perfect-clawfree} and~\ref{thm:lifts-box-perfect}.

To compute the spectrum, we apply the general spectral result for symmetric lifts of $d$-regular graphs. The complete graph $K_n$ is $(n-1)$-regular, with adjacency spectrum consisting of $\lambda_1 = n-1$ (with multiplicity $1$) and $\lambda_i = -1$ (with multiplicity $n-1$). From Theorem~\ref{thm:regular-spectrum}, the eigenvalues of the lift are:
\[
\{ d + \lambda_i - 2,\; d - \lambda_i - 2 \} = \{ n - 3 \pm \lambda_i \},
\]
for each eigenvalue $\lambda_i$ of $K_n$. Substituting, we obtain:
\[
\{ n - 3 \pm (n - 1) \} = \{ 2n - 4,\; -2 \}, \quad \text{and} \quad \{ n - 3 \pm (-1) \} = \{ n - 4,\; n - 2 \}.
\]
Thus, the full spectrum is:
\[
\{ 2(n-2),\; n - 2,\; n - 4,\; -2 \},
\]
with multiplicities $1$, $n-1$, $n-1$, and $n(n-1)-2n+1 = n^2-3n+1$ respectively.

Finally, the bipartite graph $B_{\HL'}(K_n)$ is simply the graph on $V'\sqcup V''$ with
$u'\sim v''$ iff $u\ne v$. Equivalently, it is $K_{n,n}$ with the perfect matching
$\{u'u'':u\in[n]\}$ removed, hence it is $(n-1)$-regular and connected.
Viewing $V'$ as points and $V''$ as blocks, each block $v''$ is incident with all points
except $v'$, so each block has size $n-1$, and any two distinct points $u',w'$ lie
together in exactly $n-2$ blocks (all blocks except $u''$ and $w''$). Thus $B_{\HL'}(K_n)$
is the incidence graph of the symmetric $2$-design $2\text{--}(n,n-1,n-2)$.
\end{proof}

\section{A Further Regular Family: Paley Graphs}\label{sec:paley}

\subsection{Symmetric lifts of Paley graphs}\label{subsec:paley-theory}

We record one further regular family illustrating the general spectral formula, namely the
symmetric lifts of Paley graphs. Let \(G_p\) denote the Paley graph on \(\mathbb{F}_p\), where
\(p\equiv 1 \pmod 4\) is prime: two vertices \(x,y\in\mathbb{F}_p\) are adjacent iff \(x-y\) is a
nonzero quadratic residue modulo \(p\). Then \(G_p\) is undirected and \(d\)-regular with
\(d=\frac{p-1}{2}\), and it has adjacency spectrum
\[
\Spec(G_p)=\left\{\, d,\ \left(\frac{-1+\sqrt p}{2}\right)^{\times d},\ \left(\frac{-1-\sqrt p}{2}\right)^{\times d}\right\};
\]
see, for example,~\cite{godsil2001algebraic,paley1933,brouwer2011spectra}.

Applying the symmetric lift \( \HL'_2(G_p) \) gives a graph with
\[
|V(\HL'_2(G_p))| \;=\; 2|E(G_p)| \;=\; p d \;=\; \frac{p(p-1)}{2},
\qquad
\deg(\HL'_2(G_p)) \;=\; 2(d-1) \;=\; p-3,
\]
and its adjacency spectrum is obtained directly from Theorem~\ref{thm:regular-spectrum}. This family is
included as an explicit regular example in which all numerical parameters of the lift can be written down
in closed form.

\begin{theorem}[Spectrum of \(\HL'_2(G_p)\)]\label{thm:paley-lift-spectrum}
Let \(G_p\) be the Paley graph on \(p\equiv 1\pmod 4\) vertices, with degree \(d=\frac{p-1}{2}\).
Then
\[
\Spec\!\big(\HL'_2(G_p)\big)
=
\Big(\ \{\, d-2+\lambda \,:\, \lambda\in\Spec(G_p)\,\}\ \cup\ \{\, d-2-\lambda \,:\, \lambda\in\Spec(G_p)\,\}\ \Big)
\ \cup\ \{-2\}^{\times\,p(d-2)}.
\]
\end{theorem}

\begin{proof}
Apply Theorem~\ref{thm:regular-spectrum} with \(G=G_p\), \(n=|V(G_p)|=p\), and \(d=(p-1)/2\).
The \(-2\) multiplicity is \(n(d-2)=p(d-2)\).
\end{proof}

In particular, writing the two nontrivial Paley eigenvalues as
\(\lambda_{\pm}=\frac{-1\pm \sqrt p}{2}\), the second-largest adjacency eigenvalue of the lift is
\[
\lambda_2\!\big(\HL'_2(G_p)\big)
=
\max\big\{\, d-2+\lambda_{+},\ d-2-\lambda_{-}\,\big\}
=
d-2+\lambda_{+}
=
\frac{p-5+\sqrt p}{2}.
\]
Thus the adjacency spectral gap of the \((p-3)\)-regular lift is
\[
\Gap\!\big(\HL'_2(G_p)\big)
=
(p-3)-\lambda_2\!\big(\HL'_2(G_p)\big)
=
\frac{p-1-\sqrt p}{2}.
\]
This gives an explicit pseudorandom regular family illustrating the general spectral formula together with the
clique-number statement \(\omega(\HL'_2(G_p))=(p-1)/2\). In particular, the Paley lifts furnish a concrete
infinite family of regular perfect graphs with explicitly controlled adjacency spectrum and spectral gap, making
them a natural first model for pseudorandom behavior within the symmetric-lift construction.

\section{Iteration of the Symmetric Lift}

We record one further structural feature of the symmetric lift: the construction may be iterated.
Starting from a graph \(G_0=G\), define
\[
G_{k+1}:=\HL'_2(G_k)\qquad (k\ge 0).
\]
In the regular case, both the degree sequence and the spectrum of this tower admit explicit recursion.

\subsection{Regular recursion}

Suppose \(G_0\) is \(d_0\)-regular, and write \(d_k\) for the degree of \(G_k\).
Since the symmetric lift of a \(d_k\)-regular graph is \((2d_k-2)\)-regular, we obtain
\[
d_{k+1}=2d_k-2.
\]
Equivalently, if \(d_0=d\), then
\[
d_k = 2^k(d-2)+2.
\]

If \(n_k:=|V(G_k)|\), then
\[
n_{k+1}=2|E(G_k)| = n_k d_k.
\]
Thus the number of vertices grows multiplicatively according to the degree sequence of the tower.

\subsection{Spectral recursion in the regular case}

Let \(G_k\) be \(d_k\)-regular with adjacency eigenvalues
\[
\lambda_1^{(k)}=d_k,\ \lambda_2^{(k)},\dots,\lambda_{n_k}^{(k)}.
\]
By Theorem~\ref{thm:regular-spectrum}, the adjacency spectrum of \(G_{k+1}=\HL'_2(G_k)\) is
\[
\Spec(G_{k+1})
=
\{\,d_k-2+\lambda_i^{(k)}:1\le i\le n_k\,\}
\cup
\{\,d_k-2-\lambda_i^{(k)}:1\le i\le n_k\,\}
\cup
\{-2\}^{\,n_k(d_k-2)}.
\]

Thus each non-flat eigenvalue at stage \(k\) produces two affine descendants at stage \(k+1\),
while a large flat band at \(-2\) is created with multiplicity \(n_k(d_k-2)\).

\section{Conclusion}

We have studied the canonical doubled edge-stage lift
\[
\HL'_2(G)=L(G\otimes K_2),
\]
viewed as a structural enlargement of the line graph together with a natural involution
\((u,v)\leftrightarrow (v,u)\). The central structural fact is that the quotient by this involution
recovers the ordinary line graph:
\[
\HL'_2(G)/\iota \cong L(G).
\]
Thus the symmetric lift provides a canonical doubled representation of edge-adjacency in which the
line graph appears as the symmetric sector.

This doubled representation admits a natural sector decomposition. The symmetric sector recovers
\(L(G)\), while the antisymmetric sector yields the signed graph \(\mathcal A(G)\), giving a refined edge-space
object that records whether overlaps are same-end or crossed. Thus the symmetric lift naturally combines a
canonical replacement construction with a nontrivial decomposition into symmetric and antisymmetric sectors.

We have also recorded explicit consequences in the regular case, including closed spectral formulas and
concrete families such as the complete-graph lifts and the Paley lifts. These examples show that the
doubled edge-stage construction is both structurally rigid and computationally accessible, and that it
interacts naturally with dense deterministic families as well as highly structured pseudorandom ones.

The Paley family is especially useful here: it gives an explicit infinite family of regular perfect graphs
obtained canonically from pseudorandom input, with degree, spectrum, and spectral gap all available in
closed form. In this sense, the Paley lifts provide a first concrete source of what one may reasonably view
as pseudorandom or expander-type perfect graphs arising from the present construction.

More broadly, the same mechanism suggests a natural program for random regular input families. Whenever the
base graph is a regular graph with strong spectral expansion, the symmetric lift produces a regular perfect
graph with explicitly controlled second eigenvalue and spectral gap; see
Theorem~\ref{thm:lift-second-eigenvalue} and Corollary~\ref{cor:lift-ramanujan-type}. We do not pursue that
probabilistic direction here, but it suggests that random regular graphs may yield a corresponding source of
canonical perfect graphs with strong inherited spectral expansion.

The antisymmetric sector introduced here is developed intrinsically in separate work on the antisymmetric
line graph. The present paper isolates the symmetric lift itself as the canonical doubled edge-stage object:
it contains the ordinary line graph as quotient, carries a natural involution, and supports the sector
decomposition developed above.

\section*{Acknowledgements}

The author would like to thank Gaurav Bhatnagar, Amitabha Bagchi, and Nishad Kothari for discussions and insights related to this work.

\vspace{2em}
\noindent\textbf{Author address:} \\
Hartosh Singh Bal \\
The Caravan, Jhandewalan Extn., New Delhi 110055, India \\
\texttt{hartoshbal@gmail.com}

\end{document}